\documentclass[12pt]{article}

\usepackage[a4paper,margin=2.5cm]{geometry}

\usepackage{mathrsfs}
\usepackage{latexsym}
\usepackage{amsfonts}
\usepackage{amsmath}
\usepackage{amsthm}
\usepackage{amssymb}
\usepackage[english]{babel}
\input{colordvi}
\usepackage{color}
\usepackage{multirow}
\usepackage{lineno}
\usepackage{array}
\usepackage{multirow}
\usepackage[table]{xcolor}
\usepackage{verbatim}

\usepackage{subfigure}
\usepackage{arydshln} 
\usepackage{graphicx}
\usepackage{pgf}
\usepackage{bm}
\usepackage{pstricks}

\usepackage{epsfig}

\definecolor{graytable}{HTML}{D6D5D5}

\title{ {\bf On the benchmark instances for the Bin Packing with Conflicts} }

\author{Tiziano Bacci\thanks{CNR-IASI, Via dei Taurini 19, 00185 Roma, Italia -- \texttt{\{tiziano.bacci,sara.nicoloso\}@iasi.cnr.it}}\;\;and Sara Nicoloso\footnotemark[1]}

\date{\today}

\begin{document}

\maketitle

\begin{abstract}
\noindent
Many authors, mainly in the context of the Bin Packing Problem with Conflicts, used the random graph generator proposed in ``Heuristics and lower bounds for the bin packing problem with conflicts'' [M.~Gendreau, G.~Laporte, and F.~Semet, {\em Computers \& Operations Research}, 31:347--358, 2004]. In this paper we prove that the graphs generated in this way are not arbitrary but threshold ones. Computational results show that instances with threshold conflict graphs are easier to solve w.r.t.~instances with arbitrary conflict graphs.
\end{abstract} 
 
\vfil
\vfil 

\section{Introduction}
\label{sec:introduction}

In this paper we show that a popular random graph generator \cite{GLS2004}, widely used in the context of Bin Packing Problem with Conflicts, generates very special graphs namely threshold graphs and not arbitrary ones as claimed by Gendreau et al.~\cite{GLS2004}, nor arbitrary interval ones as claimed by Sadykov and Vanderbeck \cite{SV2013}.

In Section \ref{sec:threshold} we define the threshold graphs and discuss some of their peculiar properties, in Section \ref{sec:t-generator} we present the generator defined in \cite{GLS2004} showing that it produces threshold graphs, in Section \ref{sec:bppc} we analyse the effects of using this generator on instances of Bin Packing Problem with Conflicts. Concluding remarks in Section \ref{sec:conclusions}. 

\section{Threshold graphs}
\label{sec:threshold}

A graph $G=(V,E)$ is a threshold graph if there exist a real number $d$ (the threshold) and a weight $p_x$ for every vertex $x \in V$ such that $(i,j)$ is an edge iff $(p_i + p_j)/2 \le d$ (see \cite{G1980}). W.l.o.g.~from now on we assume that $p_x \in [0,1]$ $\forall x$ (as a consequence it makes sense to choose $d \in [0,1]$). 

According to this definition it follows that a vertex $i$ is connected to all the vertices $j$ such that $p_j \le 2d - p_i$. Thus, $N(h) \supseteq N(k)$ and $\deg(h) \ge \deg(k)$ if and only if $p_h \le p_k$, where $N(x)$ denotes the set of vertices adjacent to $x$ and $\deg(x) = |N(x)|$.

A threshold graph has many peculiar properties as it is at the same time an interval graph, a co-interval graph, a cograph, a split graph, and a permutation. In addition, its complement, where $(i,j)$ is an edge iff $(p_i + p_j)/2 > d$, is a threshold graph too.

\begin{figure}[htbp] 
  \begin{center}
    \begin{pspicture}(1,0)(15,5.5)


 \multiput (1,4.228333333)(0,-0.071666667) {3} {
 \multiput (0,0)(0.071666667,0) {40} {
		 \psframe*[linecolor=graytable](0,0)(0.071666667,0.071666667)
		 }
    }

 \multiput (1,4.013333333)(0.071666667,0) {37} {
		 \psframe*[linecolor=graytable](0,0)(0.071666667,0.071666667)
		 }

 \multiput (1,3.941666667)(0.071666667,0) {36} {
		 \psframe*[linecolor=graytable](0,0)(0.071666667,0.071666667)
		 }

 \multiput (1,3.87)(0.071666667,0) {35} {
		 \psframe*[linecolor=graytable](0,0)(0.071666667,0.071666667)
		 }

 \multiput (1,3.798333333)(0.071666667,0) {33} {
		 \psframe*[linecolor=graytable](0,0)(0.071666667,0.071666667)
		 }

 \multiput (1,3.726666667)(0,-0.071666667) {3} {
 \multiput (0,0)(0.071666667,0) {31} {
		 \psframe*[linecolor=graytable](0,0)(0.071666667,0.071666667)
		 }
    }

 \multiput (1,3.511666667)(0.071666667,0) {29} {
		 \psframe*[linecolor=graytable](0,0)(0.071666667,0.071666667)
		 }

 \multiput (1,3.44)(0.071666667,0) {25} {
		 \psframe*[linecolor=graytable](0,0)(0.071666667,0.071666667)
		 }

 \multiput (1,3.368333333)(0,-0.071666667) {7} {
 \multiput (0,0)(0.071666667,0) {24} {
		 \psframe*[linecolor=graytable](0,0)(0.071666667,0.071666667)
		 }
    }

 \multiput (1,2.866666667)(0,-0.071666667) {5} {
 \multiput (0,0)(0.071666667,0) {19} {
		 \psframe*[linecolor=graytable](0,0)(0.071666667,0.071666667)
		 }
    }

 \multiput (1,2.508333333)(0.071666667,0) {12} {
		 \psframe*[linecolor=graytable](0,0)(0.071666667,0.071666667)
		 }

 \multiput (1,2.436666667)(0,-0.071666667) {4} {
 \multiput (0,0)(0.071666667,0) {11} {
		 \psframe*[linecolor=graytable](0,0)(0.071666667,0.071666667)
		 }
    }

 \multiput (1,2.15)(0,-0.071666667) {2} {
 \multiput (0,0)(0.071666667,0) {10} {
		 \psframe*[linecolor=graytable](0,0)(0.071666667,0.071666667)
		 }
    }

 \multiput (1,2.006666667)(0,-0.071666667) {2} {
 \multiput (0,0)(0.071666667,0) {7} {
		 \psframe*[linecolor=graytable](0,0)(0.071666667,0.071666667)
		 }
    }

 \multiput (1,1.863333333)(0,-0.071666667) {2} {
 \multiput (0,0)(0.071666667,0) {6} {
		 \psframe*[linecolor=graytable](0,0)(0.071666667,0.071666667)
		 }
    }

 \multiput (1,1.72)(0.071666667,0) {5} {
		 \psframe*[linecolor=graytable](0,0)(0.071666667,0.071666667)
		 }

 \multiput (1,1.648333333)(0.071666667,0) {4} {
		 \psframe*[linecolor=graytable](0,0)(0.071666667,0.071666667)
		 }

 \multiput (1,1.576666667)(0,-0.071666667) {3} {
 \multiput (0,0)(0.071666667,0) {3} {
		 \psframe*[linecolor=graytable](0,0)(0.071666667,0.071666667)
		 }
    }

 \multiput (1,4.228333333)(0.071666667,-0.071666667) {19} {
		 \psframe*[linecolor=white](0,0)(0.071666667,0.071666667)
		 }

\psline[linestyle=dashed,dash=2pt 1pt](0,2.866666666)(5.55,2.866666666)

\rput (0.5,3.048333333) {{\scriptsize $\omega(G)$}}

\psline[linestyle=dashed,dash=2pt 1pt](2.433333334,5.3)(2.433333334,-0.25)

\rput {90}(2.263333334,4.75) {{\scriptsize $\omega(G)$}}

\psline[linestyle=dashed,dash=2pt 1pt](3.86666668,5.3)(3.86666668,-0.25)
      
\rput {90}(3.69666668,4.8) {{\scriptsize 2$\omega(G)$}}



 \multiput (6.3,4.228333333)(0,-0.071666667) {2} {
 \multiput (0,0)(0.071666667,0) {58} {
		 \psframe*[linecolor=graytable](0,0)(0.071666667,0.071666667)
		 }
    }

 \multiput (6.3,4.085)(0.071666667,0) {57} {
		 \psframe*[linecolor=graytable](0,0)(0.071666667,0.071666667)
		 }

 \multiput (6.3,4.013333333)(0.071666667,0) {56} {
		 \psframe*[linecolor=graytable](0,0)(0.071666667,0.071666667)
		 }
 \multiput (6.3,3.941666667)(0.071666667,0) {55} {
		 \psframe*[linecolor=graytable](0,0)(0.071666667,0.071666667)
		 }
 \multiput (6.3,3.87)(0.071666667,0) {54} {
		 \psframe*[linecolor=graytable](0,0)(0.071666667,0.071666667)
		 }

 \multiput (6.3,3.798333333)(0,-0.071666667) {2} {
 \multiput (0,0)(0.071666667,0) {53} {
		 \psframe*[linecolor=graytable](0,0)(0.071666667,0.071666667)
		 }
    }

 \multiput (6.3,3.655)(0,-0.071666667) {2} {
 \multiput (0,0)(0.071666667,0) {52} {
		 \psframe*[linecolor=graytable](0,0)(0.071666667,0.071666667)
		 }
    }

 \multiput (6.3,3.511666667)(0,-0.071666667) {3} {
 \multiput (0,0)(0.071666667,0) {49} {
		 \psframe*[linecolor=graytable](0,0)(0.071666667,0.071666667)
		 }
    }

 \multiput (6.3,3.296666667)(0,-0.071666667) {2} {
 \multiput (0,0)(0.071666667,0) {47} {
		 \psframe*[linecolor=graytable](0,0)(0.071666667,0.071666667)
		 }
    }

 \multiput (6.3,3.153333333)(0,-0.071666667) {4} {
 \multiput (0,0)(0.071666667,0) {44} {
		 \psframe*[linecolor=graytable](0,0)(0.071666667,0.071666667)
		 }
    }

 \multiput (6.3,2.866666667)(0.071666667,0) {42} {
		 \psframe*[linecolor=graytable](0,0)(0.071666667,0.071666667)
		 }
 \multiput (6.3,2.795)(0.071666667,0) {41} {
		 \psframe*[linecolor=graytable](0,0)(0.071666667,0.071666667)
		 }
 \multiput (6.3,2.723333333)(0.071666667,0) {40} {
		 \psframe*[linecolor=graytable](0,0)(0.071666667,0.071666667)
		 }

 \multiput (6.3,2.651666667)(0,-0.071666667) {2} {
 \multiput (0,0)(0.071666667,0) {38} {
		 \psframe*[linecolor=graytable](0,0)(0.071666667,0.071666667)
		 }
    }

 \multiput (6.3,2.508333333)(0.071666667,0) {37} {
		 \psframe*[linecolor=graytable](0,0)(0.071666667,0.071666667)
		 }
 \multiput (6.3,2.436666667)(0.071666667,0) {36} {
		 \psframe*[linecolor=graytable](0,0)(0.071666667,0.071666667)
		 }
 \multiput (6.3,2.365)(0.071666667,0) {35} {
		 \psframe*[linecolor=graytable](0,0)(0.071666667,0.071666667)
		 }

 \multiput (6.3,2.293333333)(0,-0.071666667) {4} {
 \multiput (0,0)(0.071666667,0) {34} {
		 \psframe*[linecolor=graytable](0,0)(0.071666667,0.071666667)
		 }
    }

 \multiput (6.3,2.006666667)(0,-0.071666667) {3} {
 \multiput (0,0)(0.071666667,0) {31} {
		 \psframe*[linecolor=graytable](0,0)(0.071666667,0.071666667)
		 }
    }

 \multiput (6.3,1.791666667)(0.071666667,0) {27} {
		 \psframe*[linecolor=graytable](0,0)(0.071666667,0.071666667)
		 }
 \multiput (6.3,1.72)(0.071666667,0) {26} {
		 \psframe*[linecolor=graytable](0,0)(0.071666667,0.071666667)
		 }

 \multiput (6.3,1.648333333)(0.071666667,0) {25} {
		 \psframe*[linecolor=graytable](0,0)(0.071666667,0.071666667)
		 }

 \multiput (6.3,1.576666667)(0.071666667,0) {24} {
		 \psframe*[linecolor=graytable](0,0)(0.071666667,0.071666667)
		 }

 \multiput (6.3,1.505)(0,-0.071666667) {2} {
 \multiput (0,0)(0.071666667,0) {22} {
		 \psframe*[linecolor=graytable](0,0)(0.071666667,0.071666667)
		 }
    }

 \multiput (6.3,1.361666667)(0.071666667,0) {21} {
		 \psframe*[linecolor=graytable](0,0)(0.071666667,0.071666667)
		 }
 \multiput (6.3,1.29)(0.071666667,0) {20} {
		 \psframe*[linecolor=graytable](0,0)(0.071666667,0.071666667)
		 }

 \multiput (6.3,1.218333333)(0,-0.071666667) {2} {
 \multiput (0,0)(0.071666667,0) {19} {
		 \psframe*[linecolor=graytable](0,0)(0.071666667,0.071666667)
		 }
    }

 \multiput (6.3,1.075)(0,-0.071666667) {3} {
 \multiput (0,0)(0.071666667,0) {15} {
		 \psframe*[linecolor=graytable](0,0)(0.071666667,0.071666667)
		 }
    }

 \multiput (6.3,0.86)(0,-0.071666667) {2} {
 \multiput (0,0)(0.071666667,0) {13} {
		 \psframe*[linecolor=graytable](0,0)(0.071666667,0.071666667)
		 }
    }

 \multiput (6.3,0.716666667)(0,-0.071666667) {3} {
 \multiput (0,0)(0.071666667,0) {10} {
		 \psframe*[linecolor=graytable](0,0)(0.071666667,0.071666667)
		 }
    }

 \multiput (6.3,0.501666667)(0.071666667,0) {8} {
		 \psframe*[linecolor=graytable](0,0)(0.071666667,0.071666667)
		 }

 \multiput (6.3,0.43)(0.071666667,0) {6} {
		 \psframe*[linecolor=graytable](0,0)(0.071666667,0.071666667)
		 }
 \multiput (6.3,0.358333333)(0.071666667,0) {5} {
		 \psframe*[linecolor=graytable](0,0)(0.071666667,0.071666667)
		 }

 \multiput (6.3,0.286666667)(0.071666667,0) {4} {
		 \psframe*[linecolor=graytable](0,0)(0.071666667,0.071666667)
		 }

 \multiput (6.3,0.215)(0.071666667,0) {3} {
		 \psframe*[linecolor=graytable](0,0)(0.071666667,0.071666667)
		 }

 \multiput (6.3,0.143333333)(0.071666667,0) {2} {
		 \psframe*[linecolor=graytable](0,0)(0.071666667,0.071666667)
		 }

 \multiput (6.3,4.228333333)(0.071666667,-0.071666667) {31} {
		 \psframe*[linecolor=white](0,0)(0.071666667,0.071666667)
		 }

\psline[linestyle=dashed,dash=2pt 1pt](5.85,0.788333333)(10.85,0.788333333)

\rput (6.05,0.958333333) {{\scriptsize $i$}}

\psline[linestyle=dashed,dash=2pt 1pt](7.231666667,5.8)(7.231666667,-0.25)

\rput {90} (7.061666667,5.05) {{\scriptsize $last\_col(i)$}}



 \multiput (11.6,4.228333333)(0,-0.071666667) {23} {
 \multiput (0,0)(0.071666667,0) {60} {
		 \psframe*[linecolor=graytable](0,0)(0.071666667,0.071666667)
		 }
    }

 \multiput (11.6,2.58)(0,-0.071666667) {2} {
 \multiput (0,0)(0.071666667,0) {58} {
		 \psframe*[linecolor=graytable](0,0)(0.071666667,0.071666667)
		 }
    }

 \multiput (11.6,2.436666667)(0.071666667,0) {57} {
		 \psframe*[linecolor=graytable](0,0)(0.071666667,0.071666667)
		 }

 \multiput (11.6,2.365)(0,-0.071666667) {2} {
 \multiput (0,0)(0.071666667,0) {55} {
		 \psframe*[linecolor=graytable](0,0)(0.071666667,0.071666667)
		 }
    }

 \multiput (11.6,2.221666667)(0,-0.071666667) {3} {
 \multiput (0,0)(0.071666667,0) {54} {
		 \psframe*[linecolor=graytable](0,0)(0.071666667,0.071666667)
		 }
    }

 \multiput (11.6,2.006666667)(0,-0.071666667) {2} {
 \multiput (0,0)(0.071666667,0) {53} {
		 \psframe*[linecolor=graytable](0,0)(0.071666667,0.071666667)
		 }
    }

 \multiput (11.6,1.863333333)(0.071666667,0) {52} {
		 \psframe*[linecolor=graytable](0,0)(0.071666667,0.071666667)
		 }

 \multiput (11.6,1.791666667)(0,-0.071666667) {3} {
 \multiput (0,0)(0.071666667,0) {51} {
		 \psframe*[linecolor=graytable](0,0)(0.071666667,0.071666667)
		 }
    }

 \multiput (11.6,1.576666667)(0.071666667,0) {44} {
		 \psframe*[linecolor=graytable](0,0)(0.071666667,0.071666667)
		 }

 \multiput (11.6,1.505)(0.071666667,0) {41} {
		 \psframe*[linecolor=graytable](0,0)(0.071666667,0.071666667)
		 }

 \multiput (11.6,1.433333333)(0,-0.071666667) {2} {
 \multiput (0,0)(0.071666667,0) {39} {
		 \psframe*[linecolor=graytable](0,0)(0.071666667,0.071666667)
		 }
    }

 \multiput (11.6,1.29)(0,-0.071666667) {3} {
 \multiput (0,0)(0.071666667,0) {38} {
		 \psframe*[linecolor=graytable](0,0)(0.071666667,0.071666667)
		 }
    }

 \multiput (11.6,1.075)(0,-0.071666667) {7} {
 \multiput (0,0)(0.071666667,0) {37} {
		 \psframe*[linecolor=graytable](0,0)(0.071666667,0.071666667)
		 }
    }

 \multiput (11.6,0.573333333)(0.071666667,0) {34} {
		 \psframe*[linecolor=graytable](0,0)(0.071666667,0.071666667)
		 }

 \multiput (11.6,0.501666667)(0.071666667,0) {33} {
		 \psframe*[linecolor=graytable](0,0)(0.071666667,0.071666667)
		 }

 \multiput (11.6,0.43)(0.071666667,0) {31} {
		 \psframe*[linecolor=graytable](0,0)(0.071666667,0.071666667)
		 }

 \multiput (11.6,0.358333333)(0.071666667,0) {28} {
		 \psframe*[linecolor=graytable](0,0)(0.071666667,0.071666667)
		 }

 \multiput (11.6,0.286666667)(0,-0.071666667) {2} {
 \multiput (0,0)(0.071666667,0) {26} {
		 \psframe*[linecolor=graytable](0,0)(0.071666667,0.071666667)
		 }
    }

 \multiput (11.6,0.143333333)(0.071666667,0) {25} {
		 \psframe*[linecolor=graytable](0,0)(0.071666667,0.071666667)
		 }

 \multiput (11.6,0.071666667)(0,-0.071666667) {2} {
 \multiput (0,0)(0.071666667,0) {23} {
		 \psframe*[linecolor=graytable](0,0)(0.071666667,0.071666667)
		 }
    }

 \multiput (11.6,4.228333333)(0.071666667,-0.071666667) {41} {
		 \psframe*[linecolor=white](0,0)(0.071666667,0.071666667)
		 }

\psline[linestyle=dashed,dash=2pt 1pt](11.15,2.651666659)(16.15,2.651666659)

\psline[linestyle=dashed,dash=2pt 1pt](13.248333341,4.75)(13.248333341,-0.25)

\rput (11.35,2.821666659) {{\scriptsize $g$}}

\rput {90} (13.078333341,4.55) {{\scriptsize $g$}}


\rput (3.15,-1) {{\small (a)}}

\rput (8.45,-1) {{\small (b)}}

\rput (13.75,-1) {{\small (c)}}

 \multiput (1,0)(5.3,0) {3} {
		 \psframe(0,0)(4.3,4.3)
		 }

   \end{pspicture}
    \end{center}
    \vspace{1cm}
  \caption{Examples of adjacency matrices of threshold graphs with $n=60$ nodes and threshold a) $d=0.2$, b) $d=0.5$, and c) $d=0.8$.}
  \label{fig:adjacency_matrix}
\end{figure}
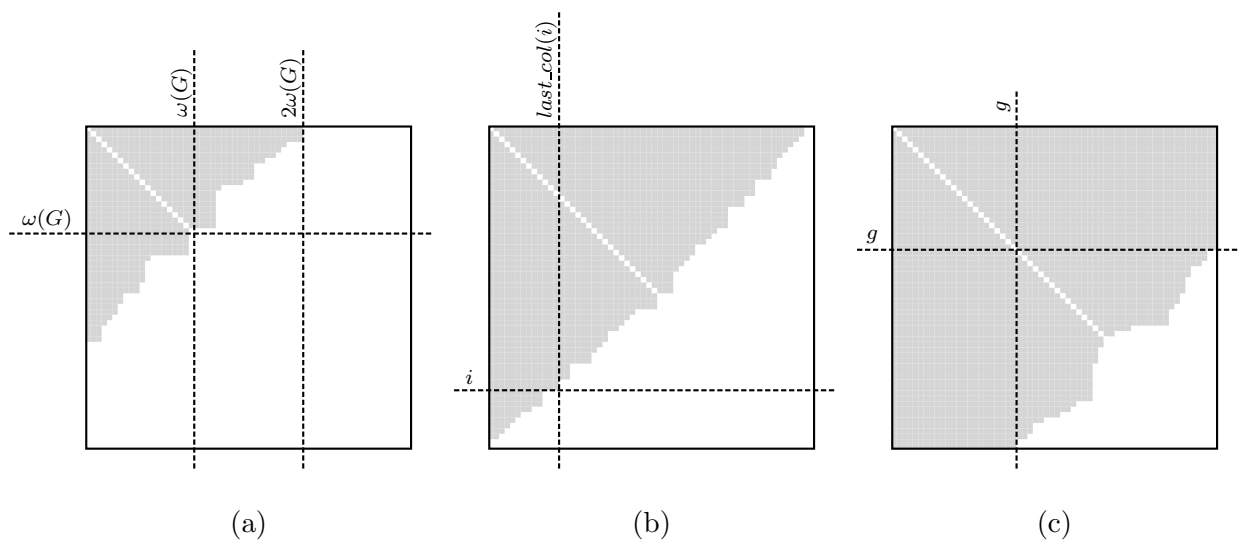

W.l.o.g.~from now on we assume that the vertices of a threshold graph $G$ are numbered in such a way that $i < j$ if and only if $\deg(i) \ge \deg(j)$. Then the $n \times n$ adjacency matrix $M = [m_{i,j}]$ of $G$ always appears as in Figure \ref{fig:adjacency_matrix}, where an entry 0 is coloured in white and an entry 1 is highlighted in grey, and $m_{i,i}=0$ for $i=1,\dots,n$. 


By what above, we observe what follows.

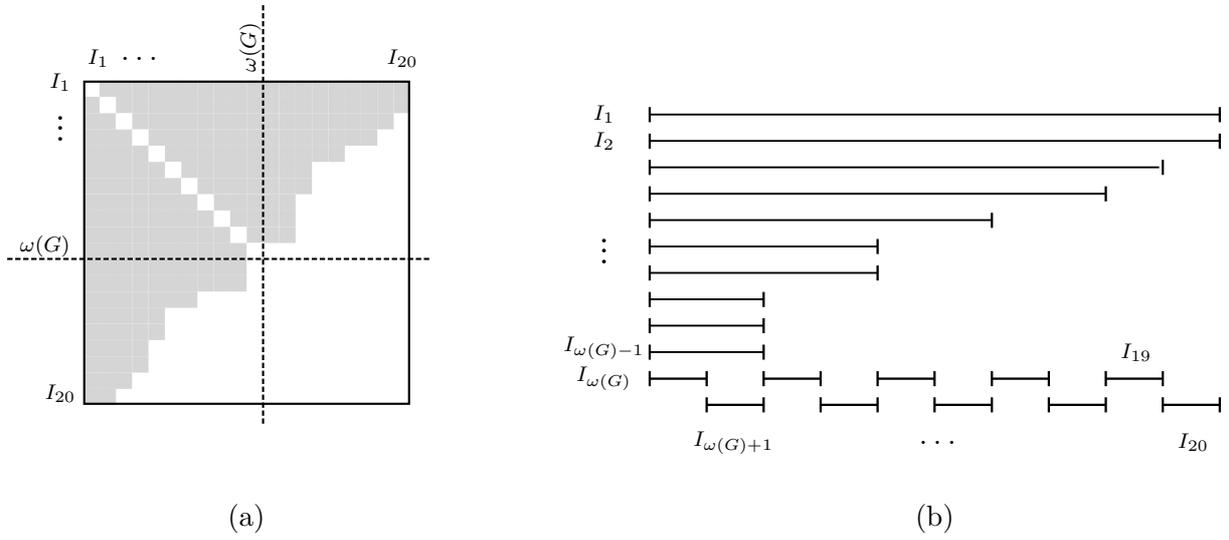
\begin{figure}[htbp] 
  \begin{center}
    \begin{pspicture}(1,0)(15,5.5)


 \multiput (1,3.87)(0,0.215) {2} {
 \multiput (0,0)(0.215,0) {20} {
		 \psframe*[linecolor=graytable](0,0)(0.215,0.215)
		 }
    }

 \multiput (1,3.655)(0.215,0) {19} {
		 \psframe*[linecolor=graytable](0,0)(0.215,0.215)
		 }

 \multiput (1,3.44)(0.215,0) {18} {
		 \psframe*[linecolor=graytable](0,0)(0.215,0.215)
		 }

 \multiput (1,3.225)(0.215,0) {16} {
		 \psframe*[linecolor=graytable](0,0)(0.215,0.215)
		 }

 \multiput (1,2.795)(0,0.215) {2} {
 \multiput (0,0)(0.215,0) {14} {
		 \psframe*[linecolor=graytable](0,0)(0.215,0.215)
		 }
    }

 \multiput (1,2.15)(0,0.215) {3} {
 \multiput (0,0)(0.215,0) {13} {
		 \psframe*[linecolor=graytable](0,0)(0.215,0.215)
		 }
    }

 \multiput (1,1.935)(0.215,0) {11} {
		 \psframe*[linecolor=graytable](0,0)(0.215,0.215)
		 }

 \multiput (1,1.505)(0,0.215) {2} {
 \multiput (0,0)(0.215,0) {10} {
		 \psframe*[linecolor=graytable](0,0)(0.215,0.215)
		 }
    }

 \multiput (1,1.29)(0.215,0) {7} {
		 \psframe*[linecolor=graytable](0,0)(0.215,0.215)
		 }

 \multiput (1,0.86)(0,0.215) {2} {
 \multiput (0,0)(0.215,0) {5} {
		 \psframe*[linecolor=graytable](0,0)(0.215,0.215)
		 }
    }

 \multiput (1,0.43)(0,0.215) {2} {
 \multiput (0,0)(0.215,0) {4} {
		 \psframe*[linecolor=graytable](0,0)(0.215,0.215)
		 }
    }

 \multiput (1,0.215)(0.215,0) {3} {
		 \psframe*[linecolor=graytable](0,0)(0.215,0.215)
		 }

 \multiput (1,0)(0.215,0) {2} {
		 \psframe*[linecolor=graytable](0,0)(0.215,0.215)
		 }

 \multiput (1,4.085)(0.215,-0.215) {11} {
		 \psframe*[linecolor=white](0,0)(0.215,0.215)
		 }

\psline[linestyle=dashed,dash=2pt 1pt](0,1.935)(5.55,1.935)

\rput (0.5,2.105) {{\scriptsize $\omega(G)$}}

\psline[linestyle=dashed,dash=2pt 1pt](3.365,5.3)(3.365,-0.25)

\rput{90} (3.195,4.75) {{\scriptsize $\omega(G)$}}

\rput (0.7,4.275) {{\scriptsize $I_1$}}

\rput (0.7,3.775) {{\large $\vdots$}}

\rput (0.7,0.15) {{\scriptsize $I_{20}$}}

\rput (1.2,4.6) {{\scriptsize $I_1$}}

\rput (1.7,4.6) {{ $\dots$}}

\rput (5.2,4.6) {{\scriptsize $I_{20}$}}





 \multiput (8.45,0.35)(1.5,0) {5} {
 \multiput (0,0)(0.75,0) {2} {
		 \psline[linewidth=0.03](0,-0.1)(0,0.1)
		 }
		 \psline[linewidth=0.03](0,0)(0.75,0)
}
 \multiput (9.2,0)(1.5,0) {5} {
 \multiput (0,0)(0.75,0) {2} {
		 \psline[linewidth=0.03](0,-0.1)(0,0.1)
		 }
		 \psline[linewidth=0.03](0,0)(0.75,0)
}

\psline[linewidth=0.025](8.45,0.7)(9.95,0.7)

 \multiput (8.45,0.7)(1.5,0) {2} {
		 \psline[linewidth=0.03](0,-0.1)(0,0.1)
		 }

\psline[linewidth=0.025](8.45,1.05)(9.95,1.05)

 \multiput (8.45,1.05)(1.5,0) {2} {
		 \psline[linewidth=0.03](0,-0.1)(0,0.1)
		 }

\psline[linewidth=0.025](8.45,1.4)(9.95,1.4)

 \multiput (8.45,1.4)(1.5,0) {2} {
		 \psline[linewidth=0.03](0,-0.1)(0,0.1)
		 }

\psline[linewidth=0.025](8.45,1.75)(11.45,1.75)

 \multiput (8.45,1.75)(3,0) {2} {
		 \psline[linewidth=0.03](0,-0.1)(0,0.1)
		 }

\psline[linewidth=0.025](8.45,2.1)(11.45,2.1)

 \multiput (8.45,2.1)(3,0) {2} {
		 \psline[linewidth=0.03](0,-0.1)(0,0.1)
		 }

\psline[linewidth=0.025](8.45,2.45)(12.95,2.45)

 \multiput (8.45,2.45)(4.5,0) {2} {
		 \psline[linewidth=0.03](0,-0.1)(0,0.1)
		 }

\psline[linewidth=0.025](8.45,2.8)(14.45,2.8)

 \multiput (8.45,2.8)(6,0) {2} {
		 \psline[linewidth=0.03](0,-0.1)(0,0.1)
		 }

\psline[linewidth=0.025](8.45,3.15)(15.15,3.15)

 \multiput (8.45,3.15)(6.75,0) {2} {
		 \psline[linewidth=0.03](0,-0.1)(0,0.1)
		 }

\psline[linewidth=0.025](8.45,3.5)(15.95,3.5)

 \multiput (8.45,3.5)(7.5,0) {2} {
		 \psline[linewidth=0.03](0,-0.1)(0,0.1)
		 }

\psline[linewidth=0.025](8.45,3.85)(15.95,3.85)

 \multiput (8.45,3.85)(7.5,0) {2} {
		 \psline[linewidth=0.03](0,-0.1)(0,0.1)
		 }

\rput (9.55,-0.5) {{\scriptsize $I_{\omega(G)+1}$}}

\rput (12.2,-0.5) {{ $\dots$}}

\rput (15.6,-0.5) {{\scriptsize $I_{20}$}}

\rput (14.85,0.7) {{\scriptsize $I_{19}$}}

\rput (7.85,0.35) {{\scriptsize $I_{\omega(G)}$}}

\rput (7.85,0.75) {{\scriptsize $I_{\omega(G)-1}$}}

\rput (7.85,2.15) {{\large $\vdots$}}

\rput (7.85,3.5) {{\scriptsize $I_2$}}

\rput (7.85,3.85) {{\scriptsize $I_1$}}

\rput (3.15,-1.5) {{\small (a)}}

\rput (12.2,-1.5) {{\small (b)}}

 \multiput (1,0)(7.2,0) {1} {
		 \psframe(0,0)(4.3,4.3)
		 }

   \end{pspicture}
    \end{center}
    \vspace{1.5cm}
  \caption{a) The adjacency matrix of threshold graph with $n=20$ nodes and b) the corresponding interval model.}
  \label{fig:Insth}
\end{figure}

\begin{enumerate}

\item \label{point1} For each row $i$, let $last\_col(i) = \max\{j : m_{i,j} = 1, j=1,\dots,n\}$ if $m_{i,1} = 1$, and $last\_col(i)=0$ if $m_{i,1} = 0$ (see Figure \ref{fig:adjacency_matrix}b); hence $last\_col(i) \ge last\_col(i+1)$.

\item Let $t =  \min \{j : m_{j,j+1} = 0, j=1,\dots,n\}$. Then the set of vertices $\{1,\dots,t\}$ induces a  maximum clique of size $\omega(G)=t$ (see Figure \ref{fig:adjacency_matrix}a). In fact, by definition, $m_{t-1,t}=1$, thus $last\_col(t-1) \ge t$ and, by Point \ref{point1},  $m_{i,j}=1$ for $i=1,\dots,t$ and $j=1,\dots,t$, $i \ne j$. 

\item  The set of vertices $\{{t},\dots,n\}$ induces a  maximum independent set of size $n-t+1$. In fact, by definition, $m_{t,t+1}=0$ and $m_{t,t-1}=1$ (as $m_{t-1,t}=1$) thus $last\_col(t) = t-1$ and $m_{i,j}=0$ for $i=t,\dots,n$ and $j=t,\dots,n$ (see Point \ref{point1}.).

\item \label{g} Let $g = \max \{h: m_{h,n}=1, h=1,\dots,n\}$ if $last\_col(1)=n$, and $g=0$ otherwise (see Figure \ref{fig:adjacency_matrix}c). Clearly, $g = last\_col(n)$. If $g \ge 1$, vertex $i$, for $i=1,\dots,g$, is connected to any other vertex.

\item \label{intervalmodel} Recalling that a threshold graph $G$ is a particular interval graph, it is always possible to derive the following family of intervals whose intersection graph is $G$: to each vertex $j= t,\dots,n$, associate the interval $I_j=(l_j,r_j)=(j-t,j-t+1)$; to each vertex \linebreak $j=1,\dots,t-1$, associate the interval $I_j=(l_j,r_j)=(0,r_{last\_col(j)})=(0,last\_col(j)-t+1)$ (we remark that $r_j \ge 1$ as $last\_col(j) \ge t\}$). See an example in Figure \ref{fig:Insth}.

\item \label{density1} The edge density $\delta = 2|E| / (n(n-1))$ of $G$ is not equal to the threshold $d$, generally speaking.

\end{enumerate}

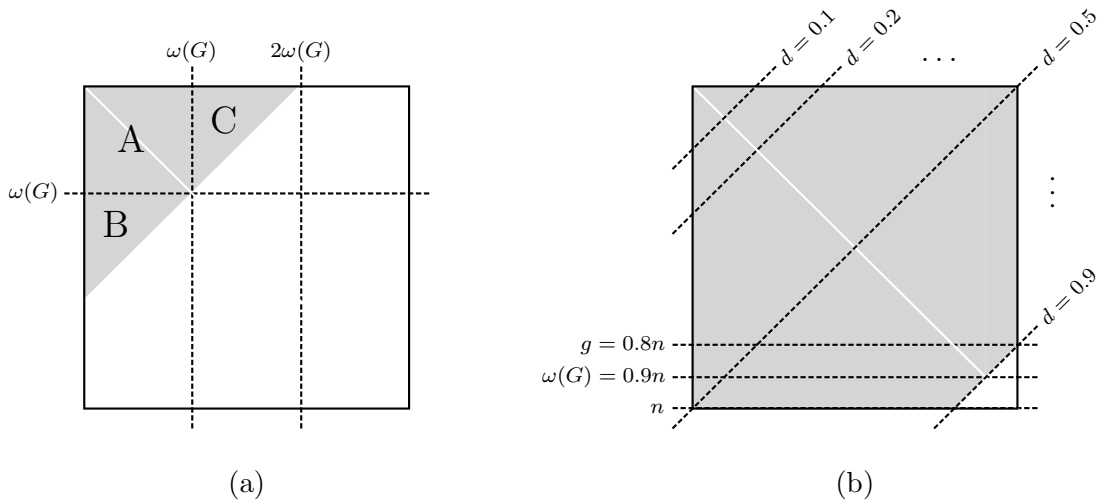
\begin{figure}[htbp] 
  \begin{center}
    \begin{pspicture}(1,0)(15,5.5)


\psframe*[linecolor=graytable](1,2.878333333)(2.421666667,4.3)

\pspolygon*[linecolor=graytable](1,1.456666666)(1,2.878333333)(2.421666667,2.878333333)

\pspolygon*[linecolor=graytable](2.421666667,2.878333333)(2.421666667,4.3)(3.843333334,4.3)

\rput (1.510833334,2.467499999) {{\large {\color{black} B }}}

\rput (2.932500002,3.839166667) {{\large {\color{black} C }}}

\psline[linestyle=dashed,dash=2pt 1pt](0.75,2.866666666)(5.55,2.866666666)

\rput (0.35,2.878333333) {{\scriptsize $\omega(G)$}}

\psline[linestyle=dashed,dash=2pt 1pt](2.433333334,4.55)(2.433333334,-0.25)

\rput (2.433333334,4.75) {{\scriptsize $\omega(G)$}}

\psline[linestyle=dashed,dash=2pt 1pt](3.86666668,4.55)(3.86666668,-0.25)
      
\rput (3.86666668,4.75) {{\scriptsize 2$\omega(G)$}}

\psline[linecolor=white](1,4.3)(5.3,0)

\rput (1.710833334,3.589166667) {{\large {\color{black} A }}}



\psframe*[linecolor=graytable](9,0.43)(12.87,4.3)

\pspolygon*[linecolor=graytable](12.44,0)(12.44,0.43)(12.87,0.43)

\pspolygon*[linecolor=graytable](12.87,0.43)(12.87,0.86)(13.3,0.86)

\psframe*[linecolor=graytable](12.87,0.86)(13.3,4.3)

\psframe*[linecolor=graytable](9,0)(12.44,0.43)

\psline[linestyle=dashed,dash=2pt 1pt](8.75,0.43)(13.55,0.43)

\rput (7.85,0.43) {{\scriptsize $\omega(G) = 0.9n$}}

\psline[linestyle=dashed,dash=2pt 1pt](8.75,0.86)(13.55,0.86)

\rput (8.1,0.86) {{\scriptsize $g = 0.8n$}}

\psline[linestyle=dashed,dash=2pt 1pt](8.75,0.02)(13.55,0.02)

\rput (8.55,0.02) {{\scriptsize $n$}}

\psline[linecolor=white](9,4.3)(13.3,0)

\psline[linestyle=dashed,dash=2pt 1pt](8.75,3.19)(10.11,4.55)

\rput {45}(10.51,4.95) {{\scriptsize $d=0.1$}}

\psline[linestyle=dashed,dash=2pt 1pt](8.75,2.33)(10.97,4.55)

\rput {45}(11.37,4.95) {{\scriptsize $d=0.2$}}

\psline[linestyle=dashed,dash=2pt 1pt](8.75,-0.25)(13.55,4.55)

\rput {45}(13.95,4.95) {{\scriptsize $d=0.5$}}

\psline[linestyle=dashed,dash=2pt 1pt](12.19,-0.25)(13.55,1.11)

\rput {45}(13.95,1.51) {{\scriptsize $d=0.9$}}

\rput (12.225,4.65) {{ \dots}}

\rput {90}(13.75,2.85) {{ \dots}}

\rput (3.15,-1) {{\small (a)}}

\rput (11.15,-1) {{\small (b)}}

 \multiput (1,0)(8,0) {2} {
		 \psframe(0,0)(4.3,4.3)
		 }

   \end{pspicture}
    \end{center}
    \vspace{1cm}
  \caption{The expected adjacency matrix when $n \rightarrow \infty$ and with threshold a) $d \le 0.5$ and b) $d \ge 0.5$}
  \label{fig:densth}
\end{figure}

\vspace{0.5cm}
\noindent
For $n \rightarrow \infty$ and $p_1,\dots,p_n$ uniformly distributed in $[0,1]$, one has:

\begin{enumerate}
\setcounter{enumi}{6
}
\item \label{omega}  $\omega(G) = t = nd$.

\item \label{density} The edge density $\delta = 2|E| / (n(n-1))$ of $G$ depends on $d$. Precisely

$$
\delta = f(d) = \left \{ \begin{array}{ll} \frac{2(nd)^2 - nd}{n(n-1)} & \mbox{for } d \le 0.5 \\ \; \\ \frac{n(n-1) - 2n^2(1-d)^2 - n(1-d)}{n(n-1)} & \mbox{for } d \ge 0.5   \end{array} \right.
$$

In fact, for $d \le 0.5$ the $2|E|$ 1's are in the area $A \cup B \cup C$ (see Figure \ref{fig:densth}a). In a similar way one can compute the number of 1's in the matrix when $d \ge 0.5$.

\item $g = 0$ when $d \le 0.5$, and $g = n(2d - 1)$ when $d \ge 0.5$ (see Figure \ref{fig:densth}b).

\end{enumerate}


\section{A random threshold graph generator}
\label{sec:t-generator}
\noindent
Gendreau et al.~\cite{GLS2004} describe the following generator: ``{\em A value $p_i$ was first assigned to each vertex $i \in V$ according to a continuous uniform distribution on $[0,1]$. Each edge $(i,j)$ of $G$ was created whenever $(p_i + p_j)/2 \le d$, where $d$ is the expected density of $G$.}''.  

This generator clearly produces threshold graphs, and the expected edge density of $G$ is not $d$ as claimed but it is the one discussed in Points \ref{density1} and \ref{density} of Section \ref{sec:threshold}. To get a threshold graph with expected edge density $\delta$ one has to set 
 
$$
d = \left \{ \begin{array}{ll} \frac{1 + \sqrt{1 + 8 n (n - 1) \delta}}{4n} & \mbox{for } \delta \le 0.5 \\ \; \\ 1 + \frac{1 - \sqrt{1 + 8 n (n - 1) (1 -\delta)}}{4n} & \mbox{for } \delta \ge 0.5   \end{array} \right.
$$
 
\noindent 
Already for $n \ge 100$ these values can be approximated to $d = \sqrt{\delta / 2}$ and \linebreak $d = 1 - \sqrt{(1 - \delta) / 2}$, respectively.

The generator by Gendreau et al.~\cite{GLS2004} has been improperly used to generate arbitrary graphs \cite{BW2005, BP2016, CFVO2015, CKHT2011, CFM2017, ELGN2011, GI2016, ThJoncour2010, JMSSV2010, JOK2015, ThKhanafer2010, KCHT2012, KCT2010, KCT2012, MG2009, MR2011, ThMuritiba2010, MIMT2010, SV2013, YLW2014}. In particular, Muritiba et al.~\cite{MIMT2010} made publicly available ``benchmark'' instances generated in this way (see http://or.dei.unibo.it/library/bin-packing-problem-conflicts) and used by many authors \cite{BP2016, CFVO2015, CKHT2011, CFM2017, ELGN2011, GI2016, ThJoncour2010, JMSSV2010, ThKhanafer2010, KCHT2012, KCT2010, ThMuritiba2010, SV2013, YLW2014}.

Most of the authors using the generator by Gendreau et al.~\cite{GLS2004} claim that they group the graphs of their test bed by edge densities, but actually they group the graphs by threshold values. Our analysis on the instances by Muritiba et al.~\cite{MIMT2010} shows that the relation between the threshold $d$ and the corresponding edge density $\delta$ is the following. 

\begin{center}
\begin{tabular}{c|c|c|c|c|c|c|c|c|c|c|c|}
\cline{2-12}
$d$     & 0 & 0.1  & 0.2  & 0.3  & 0.4  & 0.5 & 0.6  & 0.7  & 0.8  & 0.9 \\ \cline{2-12} 
$\delta$ & 0 & 0.02 & 0.08 & 0.18 & 0.32 & 0.5 & 0.68 & 0.82 & 0.92 & 0.98 \\ \cline{2-12} 
\end{tabular}
\end{center}

\noindent
We remark that the values of $\delta$ coincide with those which can be computed by the formula of Point \ref{density} in Section \ref{sec:threshold}.

\section{Bin Packing Problem with Conflicts \\ on threshold graphs}
\label{sec:bppc}

The Bin Packing Problem with Conflicts ($BPPC$), first introduced in a scheduling context by Jansen and \"Ohring \cite{JO1997}, is defined as follows: given a graph $G=(V,E)$, a nonnegative integer weight $w_i$ for each vertex $i \in V$, and a nonnegative integer $B$, find a partition of $V$ into $k$ subsets $V_1,\dots,V_k$, such that the sum of the weights of the vertices assigned to same subset is less than or equal to $B$, two vertices connected by an edge do not belong to the same subset, and $k$ is minimum. 

The minimum value of $k$ will be denoted $k_{BPPC}$. The graph $G=(V,E)$ is called {\em conflict graph} and two vertices connected by an edge are said to be {\em in conflict}. 

$BPPC$ generalizes two well known combinatorial optimization problems, the Bin Packing Problem and the Vertex Coloring Problem. In fact, $BPPC$ reduces to Bin Packing when the edge set $E$ of the graph $G$ is empty, and it reduces to Vertex Coloring when $B \ge \sum_{i \in V} w_i$ or when $G$ is complete. Observe that Vertex Coloring is solvable in linear time on threshold graphs, nevertheless $BPPC$ with a threshold conflict graph is $NP$-hard because Bin Packing is \cite{GJ1979}. 

Since threshold graphs are a subclass of interval graphs, which are in their turn a subclass of arbitrary graphs, we expect that $BPPC$ on threshold graphs is the easiest to solve. To prove our claim we conducted some computational experiments. 

By $X(n,\delta)$ with $X \in \{T,I,A\}$ we denote a set of ten instances with $n$ vertices, bound $B = 150$, and conflict graph which is a threshold graph if $X = T$, an interval graph if $X = I$, an arbitrary graph if $X = A$, with expected edge density $\delta \in \{0.02,0.08,0.18,0.32,$ $0.5,0.68,0.82,0.92,0.98\}$ (the same densities of the instances used in \cite{MIMT2010}). In particular, we choose $n \in \{250,1000\}$. 

The $T(250,\delta)$ and $T(1000,\delta)$ instances are exactly those in the classes 2 and 4 by Muritiba et al.~\cite{MIMT2010}, respectively. Precisely, given $n$, the weight of the $i$-th vertex of the $k$-th instance of $T(n,\delta)$ is the same for all $\delta$. Totally we consider 180 out of 800 of the instances by Muritiba et al.

As for the $I(n,\delta)$, the weight of the $i$-th vertex of the $k$-th instance is exactly the weight of the $i$-th vertex of the $k$-th instance of $T(n,\delta)$, and the arbitrary interval conflict graphs have been generated according to the interval graph generator by Bacci and Nicoloso \cite{BN2017IG}\footnote{The generator in \cite{BN2017IG} is not able to produce interval graphs with $n = 1000$ and edge density $\delta = 0.98$; in the corresponding cell of Table \ref{tt,ti,ta} of the present paper the average edge density of the ten instances is 0.96.}.

As for the $A(n,\delta)$, the weight of the $i$-th vertex of the $k$-th instance is exactly the weight of the $i$-th vertex of the $k$-th instance of $T(n,\delta)$, and the arbitrary conflict graphs have been generated as in \cite{SV2013}: ``{\em We began with the empty graph. We iteratively selected an item pair $(i,j)$ at random (with uniform distribution); then edge $(i,j)$ was added to the graph if it was not already defined. The procedure was interrupted as soon as the desired graph density was reached.}''.

We solve to optimality the $T(n,\delta)$, $I(n,\delta)$, and $A(n,\delta)$ instances for all $n$ and $\delta$ by means of the Vector Packing Solver 3.1.2 (VPS for short) by  Brand\~ao and Pedroso \cite{BP2016}, available at http://vpsolver.dcc.fc.up.pt/. This method is based on an arc-flow formulation with side constraints and builds very strong integer programming models that can be given in input to any state-of-the-art mixed integer programming solver (we used Cplex 12.6 on an Intel Core i7-3632QM 2.20GHz $\times$ 8 with 16 GB RAM under a Linux operating system). Actually, the arc-flow formulation is derived from a suitable graph which is preliminarily generated and whose size increases rapidly with $B$. We remark that the algorithm is applied to many classical combinatorial problems and, in particular, all the 800 instances by Muritiba et al.~\cite{MIMT2010} are solved to optimality within 50 minutes and with an average runtime of two minutes. In our analysis we set a time limit of 600 seconds for each instance. 

The computational results are summarized in Table \ref{tt,ti,ta}, where rows are indexed by $\delta$, and columns by the type of the conflict graph. In the ``Opt'' columns we report the number of instances, out of ten, solved to optimality within the time limit, and in the ``Time'' columns the time in seconds required to solve one instance, averaged over the solved instances, only.

The results in the table show that threshold instances $T$ are easier w.r.t.~instances with interval conflict graphs, and the latter are easier than instances with arbitrary conflict graphs, confirming our claim. 

We remark that, as far as we know, no tests on instances of $BPPC$ with arbitrary interval conflict graphs were performed in the literature. Sadykov and Vanderbeck \cite{SV2013} observe that the conflict graphs of the benchmark instances by Muritiba et al.~\cite{MIMT2010} are  interval graphs and not arbitrary graphs (actually they are not arbitrary interval ones). Nevertheless, to our knowledge, Sadykov and Vanderbeck \cite{SV2013} are the only ones who test their algorithm on instances with arbitrary conflict graphs.

\begin{table}[]
\centering

{\footnotesize

\begin{tabular}{cc|cc|cc|cc||cc|cc|cc|}
\cline{3-14}
                                         &      & \multicolumn{6}{c||}{$n = 250$} & \multicolumn{6}{c|}{$n = 1000$} \\ \cline{3-14}
\cline{3-14}
                                         &      & \multicolumn{2}{c|}{$T$}              & \multicolumn{2}{c|}{$I$}                & \multicolumn{2}{c||}{$A$} & \multicolumn{2}{c|}{$T$}              & \multicolumn{2}{c|}{$I$}                & \multicolumn{2}{c|}{$A$}              \\ \cline{3-14} 
                                         &      & \multicolumn{1}{c}{Opt} & Time    & \multicolumn{1}{c}{Opt} & Time      & \multicolumn{1}{c}{Opt} & Time & \multicolumn{1}{c}{Opt} & Time    & \multicolumn{1}{c}{Opt} & Time      & \multicolumn{1}{c}{Opt} & Time     \\ \hline
\multicolumn{1}{|c|}{\multirow{9}{*}{$\delta$}} & 0.02 & 10                         & 1.28 & 1                          & 138.28  & 1                          & 206.62     & 10                         & 76.37 & 0                          & -  & 0                          & -      \\ \cline{2-14} 
\multicolumn{1}{|c|}{}                   & 0.08 & 10                         & 2.75  & 0                          & -         & 0                          & -   & 10                         & 292.26  & 0                          & -         & 0                          & -      \\ \cline{2-14} 
\multicolumn{1}{|c|}{}                   & 0.18 & 10                         & 3.37 & 1                          & 522.39  & 0                          & - & 10                         & 359.16 & 0                          & -  & 0                          & -  \\ \cline{2-14} 
\multicolumn{1}{|c|}{}                   & 0.32 & 10                         & 3.81 & 10                         & 201.56 & 4                         & 340.2 & 3                         & 444.77 & 0                         & - & 0                         & - \\ \cline{2-14} 
\multicolumn{1}{|c|}{}                   & 0.5  & 10                         & 1.00 & 10                         & 15.31   & 10                         & 75.16 & 10                         & 390.94 & 0                         & -   & 0                         & - \\ \cline{2-14} 
\multicolumn{1}{|c|}{}                   & 0.68 & 10                         & 0.53 & 10                         & 3.24   & 10                         & 12.43 & 10                         & 294.12 & 0                         & -   & 0                         & -  \\ \cline{2-14} 
\multicolumn{1}{|c|}{}                   & 0.82 & 10                         & 0.29  & 10                         & 2.02     & 10                         & 5.15 & 10                         & 222.57  & 5                         & 543.92     & 0                         & -  \\ \cline{2-14} 
\multicolumn{1}{|c|}{}                   & 0.92 & 10                         & 0.11 & 10                         & 1.39  & 10                         & 2.89 & 10                         & 197.60 & 10                         & 453.98  & 0                         & - \\ \cline{2-14} 
\multicolumn{1}{|c|}{}                   & 0.98 & 10                         & 0.04 & 10                         & 1.03   & 10                         & 1.92 & 10                         & 199.36 & 10                         & 366.07   & 3                         & 561.02 \\ \hline
\end{tabular}

}

\caption{Computational results on instances with threshold ($T$), interval ($I$)$^1$, and arbitrary ($A$) conflict graphs.}
\label{tt,ti,ta}
\end{table}

\section{Concluding remarks}
\label{sec:conclusions}

In this paper we show that graphs of the $BPPC$ instances considered in \cite{BW2005, BP2016, CFVO2015, CKHT2011, CFM2017, ELGN2011, GI2016, ThJoncour2010, JMSSV2010, JOK2015, ThKhanafer2010, KCHT2012, KCT2010, KCT2012, MG2009, MR2011, ThMuritiba2010, MIMT2010, SV2013, YLW2014} and generated according to Gendreau et al.~\cite{GLS2004} are threshold graphs (and not arbitrary ones), and their edge density is not the declared one. Computational evidence suggests that $BPPC$ instances  with threshold conflict graphs are easier to solve than instances with interval or arbitrary conflict graphs. 

A consequence of using the generator by Gendreau et al.~\cite{GLS2004} in the context of $BPPC$ is that, according to Point \ref{g} in Section \ref{sec:threshold}, when $d \ge 0.5$, in any optimal solution $V_i=\{i\}$ for $i=1,\dots,g$, and $V_i$ for $i \ge g+1$ can be determined by solving a smaller instance $\mathcal Q$ defined on the last $n-g$ vertices (observe that the problem becomes simpler and simpler as $d$ increases). The conflict graph of $\mathcal Q$ is a threshold graph with expected edge density 0.5, and contains a maximum clique of expected size $(n-g)/2$. So the expected value of a lower bound for $k_{BPPC}$ on the initial instance is $g + (n-g)/2$. For example, when $n = 120$ and $d = 0.9$, $g$ is expected to be $0.8n = 96$, $\mathcal Q$ has 24 vertices, on average, and $k_{BPPC} \ge 96 + 12 = 108$ (this value appears in Table 2, column LBO, Size 120, $d = 90$ of \cite{ELGN2011}).

Finally we remark that Gendreau et al.~\cite{GLS2004} claim to use ``the procedure described in Soriano and Gendreau'' \cite{SG1996}, but this is not true. In fact, this procedure  generates ``{\em edge $(i,j)$ with probability}'' $(p_i+p_j)/2$ generalizing the uniform random graph generator and outputing arbitrary graphs.


\begin{thebibliography}{10}

\bibitem{BN2017IG}
T.~Bacci and S.~Nicoloso.
\newblock A heuristic algorithm for the Bin Packing Problem with Conflicts on Interval Graphs.
\newblock 2017. arxiv:1707.00496 [math.CO].

\bibitem{BW2005}
C.~Basnet and J.~Wilson.
\newblock Heuristics for determining the number of warehouses for storing
  non-compatible products.
\newblock {\em International Transactions in Operations Research}, 12:527--538,
  2005.

\bibitem{BP2016}
F.~Brand{\~a}o and J.~P. Pedroso.
\newblock Bin packing and related problems: General arc-flow formulation with
  graph compression.
\newblock {\em Computers \& Operations Research}, 69:56--67, 2016.

\bibitem{CFVO2015}
R.~Capua, Y.~Frota, T.~Vidal, and L.~S. Ochi.
\newblock Um algoritmo heur{\`i}stico para o problema de bin packing com
  conflitos.
\newblock Manuscript, 2015.

\bibitem{CKHT2011}
F.~Clautiaux, A.~Khanafer, S.~Hanafi, and E.~Talbi.
\newblock Le probl{\`e}me bi-objectif de bin-packing avec conflits.
\newblock Manuscript, 2011.

\bibitem{CFM2017}
D.~Cornaz, F.~Furini, and E.~Malaguti.
\newblock Solving vertex coloring problems as maximum weight stable set
  problems.
\newblock {\em Discrete Applied Mathematics}, 217:151 -- 162, 2017.

\bibitem{ELGN2011}
S.~Elhedhli, L.~Li, M.~Gzara, and J.~Naoum-Sawaya.
\newblock A {B}ranch-and-{P}rice algorithm for the bin packing problem with
  conflicts.
\newblock {\em INFORMS Journal on Computing}, 23:404--415, 2011.

\bibitem{GJ1979}
M.~R. Garey and D.~S. Johnson.
\newblock {\em Computers and Intractability: A Guide to the Theory of
  {NP}-Completeness}.
\newblock W.H. Freeman \& Co, New York, 1979.

\bibitem{GLS2004}
M.~Gendreau, G.~Laporte, and F.~Semet.
\newblock Heuristics and lower bounds for the bin packing problem with
  conflicts.
\newblock {\em Computers \& Operations Research}, 31:347--358, 2004.

\bibitem{G1980}
M.~C. Golumbic.
\newblock {\em Algorithmic Graph Theory and Perfect Graphs}.
\newblock Academic Press, 1980.

\bibitem{GI2016}
T.~Gschwind and S.~Irnich.
\newblock Dual inequalities for stabilized column generation revisited.
\newblock {\em INFORMS Journal on Computing}, 28(1):175--194, 2016.

\bibitem{JO1997}
K.~Jansen and S.~Oehring.
\newblock Approximation algorithms for time constrained scheduling.
\newblock {\em Information and Computation}, 132(2):85 -- 108, 1997.

\bibitem{ThJoncour2010}
C.~Joncour.
\newblock {\em Probl{\`e}mes de placement 2D et application {\`a}
  l'ordonnancement: mod{\'e}lisation par la th{\'e}orie des graphes et
  approches de programmation math{\'e}matique}.
\newblock PhD thesis, Universit{\'e} de {B}ordeaux {I}, July 2010.
\newblock Number 4173.

\bibitem{JMSSV2010}
C.~Joncour, S.~Michel, R.~Sadykov, D.~Sverdlov, and F.~Vanderbeck.
\newblock Column generation based primal heuristics.
\newblock {\em Electronic Notes in Discrete Mathematics}, 36:695--702, 2010.

\bibitem{JOK2015}
S.~B. Jouida, A.~Ouni, and S.~Krichen.
\newblock A multi-start tabu search based algorithm for solving the warehousing
  problem with conflict.
\newblock In H.~A.~Le Thi, T.~Pham Dinh, and N.~T. Nguyen, editors, {\em
  Modelling, Computation and Optimization in Information Systems and Management
  Sciences: Proceedings of the 3rd International Conference on Modelling,
  Computation and Optimization in Information Systems and Management Sciences -
  MCO 2015 - Part II}, pages 117--128. Springer International Publishing, 2015.

\bibitem{ThKhanafer2010}
A.~Khanafer.
\newblock {\em Algorithmes pour des probl{\`e}mes de bin packing mono- et
  multi-objectif}.
\newblock PhD thesis, Universit{\'e} Lille1, 2010.
\newblock Number 40363.

\bibitem{KCHT2012}
A.~Khanafer, F.~Clautiaux, S.~Hanafi, and E.~Talbi.
\newblock The min-conflict packing problem.
\newblock {\em Computers \& Operations Research}, 39:2122--2132, 2012.

\bibitem{KCT2010}
A.~Khanafer, F.~Clautiaux, and E.~Talbi.
\newblock New lower bounds for bin packing problems with conflicts.
\newblock {\em European Journal of Operational Research}, 206(2):281--288,
  2010.

\bibitem{KCT2012}
A.~Khanafer, F.~Clautiaux, and E.~Talbi.
\newblock Tree-decomposition based heuristics for the two-dimensional bin
  packing problem with conflicts.
\newblock {\em Computers \& Operations Research}, 39:54--63, 2012.

\bibitem{MG2009}
M.~Maiza and C.~Gu{\'e}ret.
\newblock A new lower bound for bin packing problem with general conflicts
  graph.
\newblock In {\em 1st Doctoriales STIC'09}, Universit\'e de M'sila, Alg\'erie,
  2009.

\bibitem{MR2011}
M.~Maiza and M.~S. Radjef.
\newblock Heuristics for solving the bin-packing problem with conflicts.
\newblock {\em Applied Mathematical Sciences}, 5:1739 -- 1752, 2011.

\bibitem{ThMuritiba2010}
F.~A.~E. Muritiba.
\newblock {\em Algorithms and Models For Combinatorial Optimization Problems}.
\newblock PhD thesis, Alma Mater Studiorum Universit{\`a} di Bologna, 2010.

\bibitem{MIMT2010}
F.~A.~E. Muritiba, M.~Iori, E.~Malaguti, and P.~Toth.
\newblock Algorithms for the bin packing problem with conflicts.
\newblock {\em INFORMS Journal on Computing}, 22(3):401--415, 2010.

\bibitem{SV2013}
R.~Sadykov and F.~Vanderbeck.
\newblock Bin packing with conflicts: a generic {B}ranch-and-{P}rice algorithm.
\newblock {\em INFORMS Journal on Computing}, 25(2):244--255, 2013.

\bibitem{SG1996}
P.~Soriano and M.~Gendreau.
\newblock Tabu search algorithms for the maximum clique problem.
\newblock {\em Dimacs Series in Discrete Mathematics and Theoretical Computer
  Science}, 26:221--236, 1996.

\bibitem{YLW2014}
Y.~Yuan, Y.~Li, and Y.~Wang.
\newblock An improved aco algorithm for the bin packing problem with conflicts
  based on graph coloring model.
\newblock {\em 21th International Conference on Management Science \&
  Engineering, Helsinki, Finland}, 2014.

\end{thebibliography}

\end{document}